\documentstyle[a4,11pt]{article}
\def\aa{{\scriptstyle\odot}_\nu}
\newtheorem{definition}{Definition}
\newtheorem{theorem}{Theorem}
\newtheorem{proposition}{Proposition}
\newtheorem{lemma}{Lemma}
\newtheorem{corollary}{Corollary}
\ifx\mathrm\undefined\def\mathrm#1{{\rm #1}}\fi
\makeatletter
\def\hexnumber@#1{\ifcase#1 0\or 1\or 2\or 3\or 4\or 5\or 6\or 7\or 8\or
 9\or A\or B\or C\or D\or E\or F\fi}

\font\tenmsb=msbm10 scaled \magstephalf
\font\sevenmsb=msbm7 scaled \magstephalf
\font\fivemsb=msbm5 scaled \magstephalf
\newfam\msbfam
\textfont\msbfam=\tenmsb
\scriptfont\msbfam=\sevenmsb
\scriptscriptfont\msbfam=\fivemsb
\edef\msbfam@{\hexnumber@\msbfam}
\def\Bbb#1{{\fam\msbfam\relax#1}}
\font\teneufm=eufm10 scaled \magstephalf
\font\seveneufm=eufm7 scaled \magstephalf
\font\fiveeufm=eufm5 scaled \magstephalf
\newfam\eufmfam
\textfont\eufmfam=\teneufm
\scriptfont\eufmfam=\seveneufm
\scriptscriptfont\eufmfam=\fiveeufm
\def\frak#1{{\fam\eufmfam\relax#1}}

\begin{document}

\title{Nambu mechanics, $n$-ary operations and their quantization
\thanks{The work reported here benefited from scientific contacts supported
by the EU Network Contract ERBCHRXCT940701.}}
\author{Mosh\'e FLATO \and Giuseppe DITO \and Daniel STERNHEIMER \\
D\'epartement de Math\'ematiques, Universit\'e de Bourgogne\\
BP 400, F-21011 Dijon Cedex, France}
\maketitle
\noindent{\bf Abstract.}
We start with an overview of the ``generalized Hamiltonian dynamics"
introduced in 1973 by Y. Nambu, its motivations, mathematical background
and subsequent developments -- all of it on the classical level.
This includes the notion (not present in Nambu's work) of a generalization
of the Jacobi identity called Fundamental Identity. We then briefly describe
the difficulties encountered in the quantization of such $n$-ary structures,
explain their reason and present the recently obtained solution combining
deformation quantization with a ``second quantization" type of approach on
${\Bbb R}^n$. The solution is called ``Zariski quantization" because it is
based on the factorization of (real) polynomials into irreducibles.
Since we want to quantize composition laws of the determinant (Jacobian) type
and need a Leibniz rule, we need to take care also of derivatives and this
requires going one step further (Taylor developments of polynomials over
polynomials). We also discuss a (closer to the root, ``first quantized")
approach in various circumstances, especially in the case of covariant star
products (exemplified by the case of ${\frak {su}}(2)$).
Finally we address the question of equivalence and triviality of such
deformation quantizations of a new type (the deformations of algebras are
more general than those considered by Gerstenhaber).

\section{Introduction}
In 1973, after having kept this very original idea in his files for many years
since he did not arrive at what he could have considered a really
satisfactory paper, Nambu decided to publish \cite{Na} the result ``as is".
The underlying idea of this new formalism was that in statistical mechanics
the basic result is Liouville theorem, which follows from but does not
require Hamiltonian dynamics. In geometrical terms, symplectic structures
are unimodular but the latter are more general. The lowest dimensional example
to start with is ${\Bbb R}^3$, with triplets of dynamical variables and the
Poisson bracket of two functions replaced by a ``triple bracket",
the Jacobian of three functions.

The very mention of the notion of triplets (even more so in those days) is
reminiscent of quarks for which satisfactory statistics explaining e.g.
confinement is yet to be found -- and this was in the back of Nambu's mind.
But in order to even think about such a connection one needs to solve the
problem of quantization and Nambu was not able to get very far in
solving this question. In fact quantization proved to be a very difficult
question and a solution was given only last year \cite{DFST}, at least for
polynomials on ${\Bbb R}^n$, using an elaborate construction based on
arithmetic properties of polynomials and methods of second quantization.

Nambu's paper was followed by a number of papers \cite{BF} \cite{CK}
\cite{MS} which essentially ``chilled" the idea from the physical point
of view by showing that the ``classical part" of it could be cast into
the framework of classical mechanics with Dirac constraints \cite{Dir},
for which there is a usual quantization scheme. The latter is more or less
what Nambu got (and was not satisfied with). Though perfectly correct, these
papers had for unwanted result that the whole formalism remained in limbo
for another twenty years (and we are in part responsible for that).

On the mathematical and mathematical physics side people had recurrently
been interested in Jordan algebras, which Nambu had also considered in his
abovementioned attempts to  quantization and (justly, it turns out) rejected:
nonassociativity is not the solution. More recently (re)appeared a much more
general framework, operads \cite{GK} \cite{Lo} \cite{Ma}. The solution we
shall present here is affiliated with that framework, but we shall not
develop here this point which would take us too far.

We had never quite forgotten Nambu's interesting formalism and in 1992,
independently, the underlying structure of $n$-ary brackets was studied
and the analogue of the Jacobi identity for Poisson brackets (called
{\it Fundamental Identity}, in short FI) discovered \cite{FF}  \cite{Ta}.
That identity, which somehow Nambu (and others) had missed in the 70's,
together with an explicit expression of the Leibniz rule, was the basis of
the modern developments. It is truly fundamental in many respects.
In particular it selects the really interesting generalization of Poisson
brackets to $n$-tuples from other a priori possible generalizations which
give rise to essentially trivial mathematical developments (which is why
we shall not quote any paper in those directions). From that identity emerged
the notion of Nambu-Poisson manifolds \cite{Ch} \cite{CT} \cite{Ta},
a manifold endowed with a $n$-bracket satisfying the obvious properties
of skew-symmetry and Leibniz rule and the less obvious FI. Locally such
manifolds are \cite{Ga} foliated by manifolds on which the structure is
given by a functional determinant.

In this paper we shall first present (Section~2) the earlier developments,
until and including the Fundamental Identity. Then (Section~3) we shall
describe in great lines the abovementioned solution to quantization
of such $n$-ary operations based on a subtle adaptation of deformation
quantization \cite{BFFLS} \cite{FLS} \cite{We}. The basic idea is to
replace in the computation of the determinant the ordinary product of (real)
polynomials by a symmetrized Moyal star product of their decomposition into
irreducibles (what we call their Zariski product, because of the
decomposition into irreducible polynomials) -- but the full solution is
not so simple because of linearity requirements and the fact that the
totally symmetrized Moyal product of linear functions is the usual product.	
In particular (which is consistent with the Barr triviality of Harrison
cohomology \cite{GS} for algebras of polynomials) we get algebra deformations
which are no more of the type studied by Gerstenhaber which we call here
``DrG-deformations" (as a tribute to the fundamental contribution of
``Dr. G." \cite{Ge}): for the new type of deformations we have only
${\Bbb R}$-linearity instead of ${\Bbb R}[[\nu]]$--linearity
($\nu$ being the deformation parameter), and it is a {\it generalized
deformation} in the sense that the limit for $\nu = 0$ of the deformed
algebra is the original algebra.
The algebra product which we ``deform" is here, in the spirit of second
quantization, that of an algebra of ``fields", the algebra of the
semi-group generated by irreducible polynomials (polynomials over the
set of normalized irreducible polynomials). Moreover, since we need to
quantize composition laws of the functional determinant (Jacobian) type
and need a Leibniz rule, we must also take care of derivatives and this
requires going one step further (Taylor developments of polynomials over
polynomials).

When we deal with covariant star products for which the totally
symmetrized star product of elementary factors is no more trivial, some
variants can be developed, of the ``first quantized" type; we present
them here in Section~4.1, essentially in the case of ${\frak {su}}(2)$
\cite{DF}. A similar procedure can be introduced also on
${\Bbb R}^{m}$ when the ordering used in the corresponding ``Zariski"
product (e.g. of standard type) is not the same as the ordering for the
star product. A second part of the last Section is devoted to some remarks
on the important question of equivalence of such generalized deformations
(which is not straightforward) and to the related (a little simpler) question
of triviality, already tackled in \cite{DF} where the generalized
deformations introduced are not trivial relatively to the more natural
notion of ``strong triviality" (though they are trivial for the less
demanding ``weak triviality").
The physical aspect of triviality requires to address also the
problem of spectrality of observables, by which we conclude that section.

\section{Nambu mechanics and its foundations}

\subsection {Nambu's original paper (classical part)}

Having in mind the abovementioned motivation, Nambu \cite{Na} started with
the following ``Hamilton equations" on ${\Bbb R}^3$ of the form:
\begin{equation}\label{a}
\frac{d{r}}{dt}=\nabla g({r}) \wedge \nabla h({r}),
 \quad {r}=(x,y,z)\in {\Bbb R}^3,
\end{equation}
where $x$, $y$, $z$ are the dynamical variables and $g$, $h$ are two
functions of ${r}$. Liouville theorem follows directly from the identity
$\nabla\cdot(\nabla g({r}) \wedge \nabla h({r}))=0$,
which tells us that the velocity field in Eq.~(\ref{a}) is divergenceless.
{}From (\ref{a}) follows that the evolution of a function $f$ on
${\Bbb R}^3$ is given by:
\begin{equation}\label{c}
\frac{df}{dt}=\frac{\partial(f,g,h)}{\partial(x,y,z)} ,
\end{equation}
where the right-hand side is the Jacobian of the mapping
${\Bbb R}^3 \rightarrow {\Bbb R}^3$ given by
$(x,y,z) \mapsto (f,g,h) $.

For example, in this ``baby model for integrable systems", Euler equations
for the angular momentum of a rigid body are obtained when the dynamical
variables are taken to be the components of the angular  momentum vector
$ L = (L_x, L_y, L_z)$, and $g$ and $h$ are, respectively, the total
kinetic energy ${L_x^2 \over 2I_x} + {L_y^2 \over 2I_y} + {L_z^2 \over 2I_z}$
and the square of the angular momentum $L_x^2 + L_y^2 + L_z^2$.	
Other examples can be given, in particular Nahm's equations for static
${\frak {su}}(2)$ monopoles, $\dot{x}_i = x_j x_k$ ($i,j,k = 1,2,3$) in
${\frak {su}}(2)^* \sim {\Bbb R}^3$, with $h = x_1^2 -x_2^2$,
$g = x_1^2 -x_3^2$, etc.

In this framework the analogues of canonical transformations are mappings
$(x,y,z) \mapsto (x',y',z')$ with
$\frac{\partial(x',y',z')}{\partial(x,y,z)} = 1 $.
The linear ones generate $SL(3,{\Bbb R})$. The Hamilton equations
generate infinitesimal canonical transformations, and two sets of
``Hamiltonians" $h, g$ and $h', g'$ generate the same transformation if
they are related by a ``gauge transformation"
	$\frac{\partial(h',g')}{\partial(h, g)}  =  1 $.
Here the principle of least action, which states that the classical
trajectory $C_1$ is an extremal of the action functional
$ A(C_1) = \int_{C_1} (pdq - Hdt) $, is replaced  by a similar one
\cite{Ta} with a 2-dimensional cycle $C_2$ and ``action functional''
 $A(C_2) =  \int_{C_2} (xdy \wedge dz - h dg \wedge dt) $
(which bears some flavor of strings and some similitude with
the cyclic cocycles of Connes \cite{Co}).								

Expression (\ref{c}) was easily generalized to $n$ functions $f_i,
i=1, \ldots, n$.  As was already noticed by Nambu (and in apparent
contradistinction with the definition of Poisson brackets), the correct
way to proceed (see \cite{Na}, \cite{BF}) is not by a direct sum of $m$
triples, where the linear canonical transformations are $SL(3,{\Bbb R})^m$
and the pairs of ``Hamiltonians" ($h_i, g_i)$, $i = 1, \ldots, m$, are no
longer constants of motion. One should instead introduce a $n$-tuple of
functions on ${\Bbb R}^n$ with composition law given by their Jacobian,	with
linear canonical transformations $SL(n,{\Bbb R})$ and a corresponding
$(n-1)$-form which is the analogue of the Poincar\'e-Cartan integral invariant.
The Jacobian has to be interpreted as a generalized Poisson bracket:
it is skew-symmetric with respect to the $f_i$'s and a derivation of the
algebra of smooth functions on  ${\Bbb R}^n$, i.e., the Leibniz rule is
verified in each argument (e.g. $\{f_1f_2, f_3, \ldots, f_{n+1}\} =
f_1 \{f_2, \ldots, f_{n+1}\} + \{f_1, f_3, \ldots, f_{n+1}\}f_2$, etc.).
Hence there is a complete analogy with the Poisson bracket formulation
of Hamilton equations (including the fact that the components of the
$(n-1)$-tuple of ``Hamiltonians" $(f_2, \ldots f_n)$ are constants of motion)
except, at first sight, for the equivalent of the Jacobi identity which
seems to be lacking.

\subsection {Nambu's mechanics as constrained mechanics}

As we said in the Introduction, it turns out \cite {BF} that Nambu mechanics
on ${\Bbb R}^n$ can be cast into the framework of classical mechanics with
Dirac constraints \cite{Dir} on ${\Bbb R}^{2n}$. Take $n = 3$ to fix ideas.
Motivated by (\ref{a}), which can be obtained using the intuitive Hamiltonian
${\cal H} = {p}\cdot (\nabla g({q}) \times \nabla h({q})),
 \quad {q} = {r} =(x,y,z)\in {\Bbb R}^3$,
we consider the ``Nambu Lagrangian"
$L_N = h. \sum_{i=1}^3 {\dot {q_i}} {{\partial g} \over {\partial q_i}}$
for which the Euler-Lagrange equations can be written as
$ {{dg} \over {dt}} {\nabla h} = {{dh} \over {dt}} {\nabla g}$.
Then (in an open set $\Omega \in {\Bbb R}^3$ where
${\nabla h} \times {\nabla g}$ is nonzero, by vector multiplication
by ${\nabla h}$ and ${\nabla g}$) we see that $h$ and $g$ are
constants of motion, as required.

But now in addition to the equation for $\dot q$ there is the ``other half"
for $\dot p$ and thus we need \cite{BF} to introduce three Dirac constraints
 $\phi_i(q,p) = p_i - h(q)\partial_i g(q)$, $i = 1,2,3$, where
${\partial_i} = {\partial \over {\partial q_i}}$ (these are primary
constraints which give a pair of second class constraints and one first
class constraint; there are no secondary constraints). We get
(with 0 for ``Hamiltonian") $\dot {q_i} = u_i$ and
$\dot p_i = \sum_{j=1}^3 u_i\partial_i(h(\partial_jg))$ for some functions
$u_i$ to be determined. Expressing $\dot {\phi_i} = 0$ we get
$\sum_{j=1}^3 u_i(\partial_ih\partial_jg - \partial_jh \partial_ig) = 0$
hence
\begin{equation}\label{e}
\dot {q_i} = u_i = v(t,q) [{\nabla h} \times {\nabla g}]_i
 \end{equation}
with an arbitrary function $v$ and finally, with a position-dependent
time rescaling $\tau = \int_0^t v(s,q(s,q_0))ds$ we obtain
${dq \over {d\tau}} = {\nabla h} \times {\nabla g}$.
Thus for a fixed $v$ (supposed to be nonzero in $\Omega$) we have
a one-to-one correspondence between solutions of (\ref{a}) and
(\ref{e}) which becomes an identity by a time rescaling. Since
the choice of the time axis is arbitrary, the two formulations contain
the same dynamical information.
${\cal H}_T = v(t,q)\ {p} \cdot ({\nabla h} \times {\nabla g})$
is now the total Hamiltonian and it
vanishes weakly since it is proportional to the first class  constraint
$\phi = {\cal H} = [{\nabla h} \times {\nabla g}] \cdot [p - h{\nabla g}]$
(the effect of the two second class constraints is to reduce the problem
to some ${\Bbb R}^4 \subset {\Bbb R}^6$).
We can now try and quantize along the lines of constrained mechanics
as proposed by Dirac and recover the usual Heisenberg quantization already
proposed by Nambu. But such an approach is not unambiguous and it is not
a quantization of the triple bracket. As we shall see in Section~3,
there is a much more fruitful and original approach which applies
directly to all $n$-tuple brackets ($n \geq 2$).

\subsection {Fundamental Identity and Nambu-Poisson manifolds}

\subsubsection {Nambu-Poisson brackets and the Fundamental Identity}
In the usual Poisson formulation, the Jacobi identity is the
infinitesimal form of the Poisson theorem which states that the bracket
of two integrals of motion is also an integral of  motion. If we want
a similar theorem for Nambu mechanics there must be an  infinitesimal
form of it which will provide a generalization of the Jacobi
identity. Denote by $\{f,g,h\}$ the Jacobian appearing in
(\ref{c}). Let $\phi_t\colon{r}\mapsto\phi_t({r})$ be the
flow for (\ref{a}).  Then a generalization of the Poisson theorem
would imply that $\phi_t$  is a ``canonical transformation'' for the
generalized bracket:
$$
\{f_1\circ\phi_t,f_2\circ\phi_t,f_3\circ\phi_t\}=\{f_1,f_2,f_3\}
\circ\phi_t .
$$
Differentiation of this equality with respect to $t$ yields the
desired  generalization of the Jacobi identity:
\begin{eqnarray*}
\{\{g,h,f_1\},f_2,f_3\}+\{f_1,\{g,h,f_2\},f_3\}+
\{f_1,f_2,\{g,h,f_3\}\}\\ =\{g,h,\{f_1,f_2,f_3\}\}, \quad \forall
g,h,f_1,f_2,f_3\in C^\infty({\Bbb R}^3).
\end{eqnarray*}
This identity and its generalization to ${\Bbb R}^{n}$, called
Fundamental  Identity (FI), was introduced by Flato, Fr\o nsdal
\cite{FF} and Takhtajan  \cite{Ta}  as a consistency condition for
Nambu Mechanics (this consistency  condition was also formulated in
\cite{SV}) and allows a generalized Poisson theorem: the generalized
bracket of $n$ integrals of motion is an  integral of motion. It turns
out that the Jacobian on ${\Bbb R}^{n}$ satisfies the FI. We are thus
lead to the following generalization, denoting by $S_n$ the group of
permutations of the set $\{1,\ldots,n\}$ and by
$\epsilon(\sigma)$ the sign of the permutation $\sigma\in S_n$:

\begin{definition} A Nambu bracket of order $n$ on a $m$-dimensional
manifold $M$ ($2\leq n\leq m$) is defined by a $n$-linear map on
$A = C^{\infty}(M)$ taking values in $A$:
$$
\{\cdot,\ldots,\cdot\}\colon A^n\rightarrow A,
$$
such that the following relations are verified
$\forall f_0,\ldots,f_{2n-1}\in A$:
\begin{itemize}
\item[a)] Skew-symmetry
\begin{equation}\label{f}
\{f_1,\ldots,f_n\} = \epsilon(\sigma)
\{f_{\sigma_1},\ldots,f_{\sigma_n}\}, \quad \forall \sigma\in S_n;
\end{equation}
\item[b)] Leibniz rule
\begin{equation}\label{h}
\{f_0f_1,f_2,\ldots,f_n\} =
f_0\{f_1,f_2,\ldots,f_n\}+\{f_0,f_2,\ldots,f_n\}f_1;
\end{equation}
\item[c)] Fundamental Identity
\begin{eqnarray}\label{i}
&&\{f_1,\ldots,f_{n-1},\{f_{n},\ldots,f_{2n-1}\}\}\nonumber\\ &&\quad
=\{\{f_1,\ldots,f_{n-1},f_n\},f_{n+1},\ldots,f_{2n-1}\}\nonumber\\
&&\qquad
+\{f_n,\{f_1,\ldots,f_{n-1},f_{n+1}\},f_{n+2},\ldots,f_{2n-1}\}
\nonumber\\ &&\qquad +\cdots+\{f_{n},f_{n+1},\ldots,f_{2n-2},
\{f_1,\ldots,f_{n-1},f_{2n-1}\}\}.
\end{eqnarray}
\end{itemize}
\end{definition}
Properties a) and b) imply that there exists a $n$-vector field $\eta$
 on $M$ such that:
\begin{equation}\label{aa}
\{f_1,\ldots,f_n\}=\eta(df_1,\ldots,df_n), \quad \forall
f_1,\ldots,f_n\in A.
\end{equation}

\subsubsection {Nambu-Poisson manifolds and associated Nambu mechanics}
Of course the FI imposes constraints on $\eta$, analyzed in \cite{Ta}.
A $n$-vector field on $M$ is called a Nambu tensor if its associated
Nambu bracket defined by (\ref{aa}) satisfies the FI. This brings us to:

\begin{definition} A Nambu-Poisson manifold $(M,\eta)$ is a manifold $M$ on
 which is defined a Nambu tensor $\eta$.  Then $M$ is said to be
endowed with a Nambu-Poisson structure.
\end{definition}
The dynamics associated with a Nambu bracket on $M$ is specified by
$n-1$ Hamiltonians $H_1,\ldots,H_{n-1}\in A$ and the time evolution of
$f\in A$ is given by:
\begin{equation}\label{j}
\frac{df}{dt}=\{H_1,\ldots,H_{n-1},f\}.
\end{equation}
Suppose that the flow $\phi_t$ associated with (\ref{j}) exists
and let $U_t$ be the one-parameter group acting on $A$ by $f\mapsto
U_t(f)=f\circ\phi_t$. It follows from the FI that:
\begin{proposition}
The one-parameter group $U_t$ is an automorphism of the Nambu bracket
structure on $A$.
\end{proposition}
\begin{definition}
$f\in A$ is called an integral of motion for the system defined by
 (\ref{j}) if it satisfies $\{H_1,\ldots,H_{n-1},f\}=0$.
\end{definition}
It follows from the FI that a Poisson-like theorem exists for
Nambu-Poisson manifolds:
\begin{proposition}
The Nambu bracket of $n$ integrals of motion is also an integral of
motion.
\end{proposition}

For the case $n=2$, the FI is the Jacobi identity and one recovers the
usual definition of a Poisson manifold.  On ${\Bbb R}^2$, the canonical
Poisson bracket of two functions ${\cal P}(f,g)$ is simply their
Jacobian, and Nambu defined his bracket on ${\Bbb R}^n$ as a Jacobian
of $n$ functions $f_1,\ldots,f_n\in C^\infty({\Bbb R}^n)$ of $n$
variables $x_1,\ldots,x_n$:
\begin{equation}\label{j2}
\{f_1,\ldots,f_n\} = \sum_{\sigma\in S_n}\epsilon(\sigma)
\frac{\partial f_1}{\partial x_{\sigma_1}}\cdots\frac{\partial f_n}
{\partial x_{\sigma_n}},
\end{equation}
which gives the {\it canonical Nambu bracket} of order $n$ on ${\Bbb R}^n$.
Other examples of Nambu-Poisson structures have been found \cite{CT}.
One of them is a generalization of linear Poisson structures and is
given by the following Nambu bracket of order $n$ on ${\Bbb R}^{n+1}$:
\begin{equation}\label{j3}
\{f_1,\ldots,f_n\} =\sum_{\sigma\in S_{n+1}}\epsilon(\sigma)
\frac{\partial f_1}{\partial x_{\sigma_1}}\cdots\frac{\partial f_{n}}
{\partial x_{\sigma_{n}}}x_{\sigma_{n+1}}.
\end{equation}
In general any manifold endowed with a Nambu-Poisson structure of
order~$n$ is {\it locally foliated} \cite{Ga} by Nambu-Poisson manifolds
of dimension~$n$ endowed with the canonical Nambu-Poisson structure
(\ref{j2}). In particular it is shown in \cite{Ga} that any Nambu
tensor is decomposable (this fact was conjectured by Takhtajan,
eventually proved in a very elegant way by Gautheron and then discovered
by chance by Takhtajan to follow from a result going back to 1923,
reproduced in a textbook by Schouten!).

We now have a complete parallel between Poisson classical mechanics
and Nambu classical mechanics: we have a bracket (the Nambu $n$-bracket,
with $n = 2$ in the former case) satisfying the required properties
(skew symmetry, Leibniz	rule and the FI), equations of evolution
of a similar form (\ref{j}) with solutions which can be written as
$f_t = {\exp}(t\{H_1, \ldots , H_{n-1}, . \})f_0$, and an {\it action
functional} for which the classical trajectory is an extremal \cite{Ta}.
That action is the integral  $A(C_{n-1}) = \int_{C_{n-1}} \omega^{(n-1)}$
over $(n-1)$-chains $C_{(n-1)}$ in the extended phase space
$M \times {\Bbb R}$ 
of the generalized Poincar\'e-Cartan integral invariant action form
$\omega^{(n-1)} = x_1dx_2 \wedge \ldots \wedge dx_n -
H_1 dH_2 \wedge \ldots \wedge dH_{n-1} \wedge dt$. We refer to \cite{Ta}
for more details (especially for the case $n=3$).

\section {Quantization of Nambu mechanics}

\subsection {First attempts and problems raised}
As realized already by Nambu (\cite{Na}), the first problem in the intriguing
question of quantization (of such a system of equations) is how to define
a quantization.	A simplistic answer is to say that, since the classical
system can be treated as a (Dirac) constrained Hamiltonian system, it is
enough to quantize the latter. Even within this framework a rigorous
solution (in particular due to the presence of a first class constraint
which turns the physical manifold into a Poisson, not symplectic, manifold)
is not so easy to obtain. But there is more: how can we be sure that the
classical equivalence can be carried over to the quantum case? Quantum
anomalies pop up like mushrooms nowadays, so one must attack directly the
problem of quantization of the original classical system.

Nambu \cite{Na}, assuming only (\ref{f}) and (\ref{h}), could not find
a satisfactory solution. Assuming the ``usual framework" (operator algebras)
to quantize (\ref{j}) he arrives at the ``disappointing" (albeit not so
surprising in view of the classical equivalence with a constrained system)
result that the quantum analogue of (\ref{j}) is equivalent to a Heisenberg
equation. Indeed, skew-symmetry~(\ref{f}) gives what is called in \cite{Ta}
($n=3$ case) the Nambu--Heisenberg commutation relation:
\begin{equation}\label{NH}
[A_1,A_2,A_3]\equiv \sum_{\sigma\in S_3}
 \epsilon(\sigma)A_{\sigma_1}A_{\sigma_2}A_{\sigma_3}=cI,
\end{equation}
where $c$ is a constant and $I$ is the unit operator. Representations
could be found (a more complete treatment is given in \cite{Ta}), but
already the Leibniz rule (\ref{h}) brings back to the Heisenberg case.
Using associative algebras to find realizations of the Poisson
bracket for the canonical triplet leads again nowhere. He thus naturally
turns to the use of nonassociative algebra. The Cayley--Dickson algebra
(an algebra with alternating associator and 8 generators) is similarly
rejected. So Nambu turns to Jordan algebras with the same conclusion, which is
no surprise since all but one of them (the exceptional algebra, $3 \times 3$
matrices with Cayley coefficients) arise from associative algebras and thus
bring back to Heisenberg formalism.	Takhtajan \cite{Ta}, with the additional
constraint of the FI, could certainly do no better. His main positive
contribution there (for the quantization problem) was to find an explicit
realization of the Nambu--Heisenberg relation (the case of general $n$
was also discussed).

\subsection {Deformation Quantization}
Deformation quantization was introduced \cite{FLS} a few years after
Nambu's paper. It is well adapted to constrained systems but (for the same
reasons as in the previous subsection) using it in this way for Nambu
mechanics will not give anything new. A more radical approach is needed
if one wants to get a solution expressing the specificity of Nambu
mechanics.

As we have mentioned in Section~2.3.2, the Jacobian is the typical case
of Nambu bracket. It is thus natural to try and replace in
(\ref{j2}) usual products of observables by star products.
For the sake of self-completeness we shall present here a very
brief overview of deformation quantization (see also \cite{FS} \cite{We}).

Let $M$ be a Poisson manifold, a differentiable manifold endowed with a
(possibly degenerate) 2-tensor $\tau$ satisfying $[\tau, \tau] = 0$
in the sense of Schouten-Nijenhuis brackets. The latter ensures that the
corresponding Poisson bracket ${\cal P}(f,g) = \tau(df \wedge dg)$
of two functions $f,g\in A = C^\infty(M)$ is a Lie algebra law.
We denote by $A[[\nu]]$ the space of formal power series in
the parameter $\nu$ with coefficients in $A$.  A {\it star product}
$\ast_\nu$ on $M$ is an associative (non-Abelian)
deformation of the usual product of the algebra $A$ defined on $A$ by:
\begin{equation} \label{st}
f\ast_\nu g = \sum_{r\geq0}\nu^r C_r(f,g), \quad \forall f,g\in A,
\end{equation}
where $C_0(f,g)=fg$, $f,g\in A$, and $C_r\colon A\times A \rightarrow
A$ ($r\geq 1$) are bidifferential operators (bipseudodifferential
operators can sometimes be considered) on $A$ satisfying:
\begin{itemize}
\item[a)] $\ast_\nu$ extends by linearity in $\nu$ to $A[[\nu]]$ and the
associativity condition $(f\ast_\nu g)\ast_\nu h = f\ast_\nu (g\ast_\nu h)$
is satisfied for all $f,g,h \in A[[\nu]]$.
\item[b)] $C_1(f,g)-C_1(g,f)=2{\cal P}(f,g)$, $f,g\in A$.
\end{itemize}
(\ref{st}) and condition a) express that we have an associative
algebra deformation in the sense of Gerstenhaber \cite{Ge} (in short,
a {\it DrG-deformation}) while condition b) ensures that the corresponding
commutator $[f,g]_{\ast_\nu}\equiv (f\ast_\nu g -g \ast_\nu f)/2\nu$
is a DrG-deformation of the Lie algebra $(A, {\cal P})$.

{\it Equivalence} of two DrG-deformations is defined in the associative
case (two star products $\ast_\nu$ and $\ast_\nu^\prime$) by the existence
of a formal series of (differential) operators $T=\sum_{r\geq0}\nu^r T_r$
with $T_0=Id$ such that
$Tf\ast Tg = T(f\ast^\prime g),\quad f,g\in A[[\nu]]$.
By equivalence one may consider only star products ``vanishing on
constants" (i.e. $C_r(f,c)=C_r(c,f)=0$, $r\geq 1$, $c\in\Bbb R$, $f\in A$,
since by \cite{GS} a deformation of a unital algebra is unital and equivalent
to a deformation with the same unit) and assume $C_1 = {\cal P}$.	
A deformation is said trivial if it is equivalent to the original product
(resp. bracket, for Lie algebra deformations); equivalence classes of
DrG-deformations are classified, at each order in $\nu$, by the second
Hochschild (resp. Chevalley) cohomology spaces $H^2(A)$
(resp. ${\tilde H}^2(A)$) of $A$ valued in itself with the natural action
(product, resp. adjoint).

When $\tau$ is everywhere nondegenerate, $M$ is {\it symplectic}
(we denote by $\omega$ the symplectic 2-form on $M$, inverse of $\tau$).
In this case (and similarly for regular Poisson manifolds; cf. e.g.
\cite{Fe} \cite{Da} \cite{We}) there always exist a star product,
${\cal P}$ is a nontrivial Hochschild 2-cocycle and
${\mathrm {dim}}{\tilde H}^2(A) = 1 + b_2$
 (where $b_2 = {\mathrm {dim}}H^2(M)$,
the second Betti number of $M$), which classifies at each order in
$\nu$ inequivalent 1-differentiable deformations \cite{FLS1} of
$(A, {\cal P})$ vanishing on constants (they are given by deformations
of the symplectic structure).  Then \cite{Fe} \cite{De} equivalence
classes of star products are in one-to-one correspondence with series
$[\omega] + \nu [H^2(M)][[\nu]]$.
In the ``flat case" $b_2 = 0$, star products are (up to equivalence)
unique and given by the Moyal product $\ast^{}_{M}$
\cite{BFFLS} which can be expressed on ${\Bbb R}^{2n}$ by the well-known
formula $f\ast^{}_{M}g=\exp(\nu{\cal P})(f,g)$. It corresponds to the Weyl
(totally symmetric) ordering of operators in quantum mechanics, the
deformation parameter being here $\nu=i\hbar/2$.

A given Hamiltonian $H\in A$ determines the time evolution of an
observable $f\in A$ by the Heisenberg equation
$\frac{df_t}{dt}=[H,f_t]_{\ast_\nu}$, the corresponding one-parameter group
of time evolution being given by the star exponential defined by:
\begin{equation}\label{exp}
\exp_{\ast}\left(\frac{tH}{2\nu}\right) \equiv \sum_{r\geq0}
\frac{1}{r!}  \left(\frac{t}{2\nu}\right)^{r}(\ast H)^r,
\end{equation}
where $(\ast H)^r=H\ast\cdots\ast H$ ($r$ factors).  Then the solution
to the Heisenberg equation above can be expressed as
$f_t = \exp_{\ast}(tH/2\nu)\ast f \ast \exp_{\ast}(-tH/2\nu)$.
In many examples (cf. e.g. \cite{BFFLS}), the star exponential is convergent
as a series in the variable $t$ in some interval ($|t|< \pi$ for the harmonic
oscillator in the Moyal case) and converges as a distribution on $M$
for fixed $t$.  Then it makes sense to consider a Fourier-Dirichlet
expansion of the star exponential:
\begin{equation}\label{fd}
\exp_{\ast}\left(\frac{tH}{2\nu}\right)(x)= \int \exp(\lambda t/2\nu)
d\mu(x;\lambda),\quad x\in M,
\end{equation}
the ``measure'' $\mu$ being interpreted as the Fourier transform (in
the distribution sense) of the star exponential in the variable $t$.
Equation (\ref{fd}) permits to define \cite{BFFLS} the spectrum of
the Hamiltonian $H$ as the support $\Lambda$ of the measure $\mu$. In
the discrete case, $\exp_{\ast}(tH/2\nu)(x)=\sum_{\lambda\in \Lambda}
\exp(\lambda t/2\nu)\pi_\lambda(x)$, $x\in M$.
The functions $\pi_\lambda$ on $M$ are interpreted as eigenstates of
$H$ associated with the eigenvalues~$\lambda$, and satisfy
$H\ast\pi_\lambda =\pi_\lambda \ast H=\lambda \pi_\lambda$,	with
$\pi_\lambda \ast \pi_{\lambda'}=\delta_{\lambda\lambda'}$ and
$\sum_{\lambda\in \Lambda} \pi_\lambda = 1$.

In the Moyal case, the Feynman path integral can be expressed
\cite{Pa} as the Fourier transform over momentum space of the
star exponential.  In field theory, where the normal star product
(which is the exponential of ``half of the Poisson bracket'' in the
variables $p\pm iq$) is relevant, the Feynman path integral is given
(up to a multiplicative factor) \cite{Di} by the star exponential.

We have limited ourselves here to some highlights (most of which will be
needed later) of deformation quantization to give the flavor of how it
provides a completely autonomous quantization scheme of a classical
Hamiltonian system. See \cite{St} for a more extensive list of relevant
references on the subject. We shall now use it as a tool for the
quantization of Nambu-Poisson structures.

\subsection {Zariski product: a generalized Abelian deformation}
Our starting point here is a simple remark: the Jacobian
of $n$ functions on ${\Bbb R}^n$ is a Nambu bracket because the usual
product of functions is Abelian, associative, distributive and
respects the Leibniz rule. This is what permits us to work out the
functional determinant and the required properties of a Nambu bracket
(including the FI) will be satisfied. Therefore, if we replace the usual
product in the Jacobian by a new product having the preceding properties,
we get a ``modified Jacobian'' which is still a Nambu bracket in the
sense of Definition 1 -- but which has ``built in" quantization if the
new product has it. In \cite{DFST} the interested reader will find
more in detail how we arrived at the solution which we shall now describe.

\subsubsection {Deformation of the usual product law}
To fix ideas, take the space $N$ of polynomials on ${\Bbb R}^3$ and denote
by $\ast$ the Moyal star product in (e.g.) the first two variables.
To get an Abelian product the natural idea is to symmetrize the Moyal
product. But if one does it brutally (i.e. ${1\over 2}(a\ast b + b\ast a),
\ a,b \in N$) one loses associativity (this is the star analogue of the
Jordan algebra formalism already studied by Nambu).

One can then try, using linearity, to bring down the problem to
monomials and then to perform total symmetrization in star monomials
(star products of the coordinates functions): symmetrization will give
back the usual product.  (This however would not necessarily be true with
different orderings or nonlinear elementary factors; we shall come back
to this point in Section~4).

But algebra tells us there is another way to decompose polynomials:
into irreducible factors, which can be of any degree. The latter is true
for polynomials in several complex variables and a fortiori in the real
domain, which is of interest for us since we want e.g. the harmonic
oscillator Hamiltonian $p^2 + q^2$ to be ``elementary". Symmetrization
adapted to the number of irreducible factors will restore associativity
of the product, but at the same time there will be ``built in" quantization
because the factors can be of any degree. That decomposition is the core
of Zariski topology on the manifold (${\Bbb R}^3$ here), hence the name
of Zariski product.

More precisely, we start with $N={\Bbb R}[x_1,x_2,x_3]$ and consider
${\cal S}(N)=\bigoplus^{\infty}_{n=1}N^{n\atop\otimes}$,
its {\it symmetric\/} tensor algebra without scalars (a kind of Fock space
over $N$). We denote by $N_1$ the semi-group of normalized polynomials,
those for which the ``maximal" monomial (defined in a natural manner:
highest total degree, then maximal with respect to the lexicographic
ordering in the variables $(x_1, x_2, x_3)$) has coefficient 1, with the
convention that $0 \in N_1$. In $N[\nu]$ we consider similarly the
semi-group (under usual product) $N^{\nu}_1$ of normalized polynomials,
those for which the non-vanishing coefficient of lowest degree in
$\nu$ belongs to $N_1$.
Now we define a map $\alpha: N^{\nu}_1 \rightarrow {\cal S}(N)$ by
$\alpha(P) = P_1\otimes\cdots\otimes P_n$ where $P_1 \cdots P_n$ is the
decomposition into irreducible factors of the coefficient of degree $0$
in $\nu$ of $P \in N^{\nu}_1$ and $\otimes$ denotes the {\it symmetric\/}
tensor product. In other words, $\alpha$ takes into account only the
classical part (which can be 0) of the polynomial in $N^{\nu}_1$, and
then replaces the usual product by the symmetric tensor product.

In order to get our ``built in" quantization we use now an {\it ``evaluation
map"} $T: {\cal S}(N) \rightarrow N[\nu]$ defined by replacing the symmetric
tensor product by a symmetrized Moyal product. In the above notations the
result of all this is
\begin{equation}\label{ta}
T (\alpha (P)) = \frac{1}{n!}\sum_{\sigma\in S_n}
P_{\sigma_1}\ast \cdots \ast P_{\sigma_n}\,
\end{equation}
and we can define an Abelian product by
\begin{equation}\label{ev}
P \times_\alpha Q=T(\alpha(P) \otimes \alpha(Q)),\quad  \forall P,Q
\in N^{\nu}_1.
\end{equation}
Obviously the product (\ref{ev}) maps $N^{\nu}_1 \times N^{\nu}_1$
into $N^{\nu}_1$ and if $P,Q \in N_1$, then $P \times_\alpha Q|_{\nu=0}=PQ$.
In this sense, $\times_\alpha$ is a generalized deformation of the usual
product on the semi-group $N_1$, but this is far from being a DrG-deformation
of the algebra $N$. In \cite{DFST} simple examples are given to make this
more explicit. We only have a deformed Abelian semi-group law, and if we plug
this law into the Jacobian then the Leibniz rule and the FI will fail.
What is lacking (and obviously so since addition and decomposition into
irreducible factors of polynomials are highly noncommutative operations)
is distributivity of the product $\times_\alpha$ with respect to addition.

In fact it follows from general results of Harrison cohomology (cf. e.g.
\cite{GS} and references therein; compare also with Gelfand's theorem on
the realization of Abelian involutive Banach algebras as algebras of
functions over the spectrum) that Abelian DrG-deformations of algebras of
polynomials (or even functions \cite{PG}) are trivial. To overcome this
difficulty and restore distributivity we use a trick inspired by the second
quantization procedure. Formally, let us look at ``functions'' on $N_1$
(e.g. formal series). Intuitively we get a deformed coproduct and the dual
of this space of ``functions'' (polynomials on polynomials) will then have
a product and a deformed product, both of which will be distributive with
respect to the vector space addition. Now the product of polynomials is again
a polynomial. So in the end we are getting some deformed product on an
algebra generated by the polynomials.

\subsubsection {The Zariski algebra and its deformations}
The ``Zariski" algebra generated by the polynomials is nothing but the
algebra ${\cal Z}_0$ of the semi-group $N_1$, i.e. the free Abelian algebra
generated by the set of normalized real irreducible polynomials
$N^{irr}_1 \subset N_1$ as building blocks.
If we denote by $Z_u$ the element of ${\cal Z}_0$ defined by $u \in N_1$,
the {\it classical Zariski product} $\bullet^z$ in ${\cal Z}_0$ is given by
$Z_u \bullet^z Z_v = Z_{uv}$, $u,v \in N_1$. Defining $Z_{cu} = cZ_u$ for
$c \in {\Bbb R}$, we can extend the above product to all $u \in N$ and get
a multiplicative injection of $N$ into ${\cal Z}_0$ which is non-additive,
i.e. $Z_{u+v}\neq Z_{u}+Z_{v}$ (the addition in ${\cal Z}_0$ is {\it not}
related to the addition in $N$). It is the algebra ${\cal Z}_0$ that we
shall now quantize.

We consider the space ${\cal Z}_\nu = {\cal Z}_0[\nu]$ of polynomials
in $\nu$ with coefficients in ${\cal Z}_0$ and inject $N_1^\nu$ into
${\cal Z}_\nu$ by	$\zeta : (\sum_{r\geq0}\nu^r u_r) \mapsto
\sum_{r\geq0}\nu^r Z_{u_r}$ with $u_0\in N_1$ and $u_i\in N, i\geq 1$.
Using this injection we can define a {\it deformed Zariski product}
$\bullet^z_\nu$ based on $\times_\alpha$, first on basis elements $Z_u$,
$u = u_1 \cdots u_m$, $u_j \in N^{irr}_1$, $j = 1, \ldots, m$,
and $Z_v$, by $Z_u \bullet^z_\nu Z_v = \zeta(u \times_\alpha v)$; then
we extend it by linearity to all of ${\cal Z}_\nu$, taking into account
the requirement that the product $\bullet^z_\nu$ annihilates the
non-zero powers of $\nu$:
\begin{equation}\label{dzp}
(\sum_{r\geq0}\nu^r A_r)\bullet^z_\nu (\sum_{s\geq0}\nu^s B_s) =
A_0\bullet^z_\nu B_0,\quad \forall A_r, B_s \in {\cal Z}_0, r,s \geq 0.
\end{equation}
The $\bullet^z_\nu$ product is obviously associative (because $\times_\alpha$
is), distributive with respect to addition in ${\cal Z}_\nu$, and Abelian.
Its limit for $\nu = 0$ is $\bullet^z$. Thus we have achieved a first step:

\begin{theorem}
The vector space ${\cal Z}_\nu$ endowed with the product $\bullet^z_\nu$
is an Abelian algebra which is a (generalized) deformation of the Abelian
algebra $({\cal Z}_0,\bullet^z)$.
\end{theorem}

The algebra ${\cal Z}_\nu$ with product $\bullet^z_\nu$ provides
thus an Abelian deformation of the algebra ${\cal Z}_0$, a fact interesting
{\it per se\/} because it gives a first example of a non-trivial Abelian
deformation, however generalized and therefore not necessarily classified by
the Harrison cohomology (defined on the sub-complex of the Hochschild complex
consisting of symmetric cochains \cite{GS}).

\subsection {Arithmetic ``second" quantization of Jacobian in ${\Bbb R}^3$}
\subsubsection {Problems with the definition of derivatives}
The next step towards our goal of deforming the Nambu bracket on ${\Bbb R}^3$
is to define the {\it derivatives} $\delta_i$, $1\leq i\leq 3$, on
${\cal Z}_0$ and extend them to ${\cal Z}_\nu$. This would allow to define
first the classical Nambu bracket on ${\cal Z}_0$, and the quantum one on
${\cal Z}_\nu$. But this is not so simple to achieve.
The straightforward definition $\delta_i Z_u=Z_{\partial_i u}$,
$\forall u\in N$, where $\partial_i$ is the usual derivative with respect
to $x^i$, does not satisfy the Leibniz rule because of the different nature
of the addition in $N$ and in ${\cal Z}_0$ (except on the diagonal, a remark
relevant for the analogue of the star exponential (\ref{exp}), which we shall
consider in Section~4.3.1).

The next natural choice is then to ``build in" the Leibniz rule by using
the	above definition of the derivatives $\delta_i$ only for irreducible
polynomials $u\in N^{irr}_1$ and postulating the Leibniz rule on the
product $v=v_1v_2\cdots v_m$ of irreducible polynomials :
$$
\delta_i\ Z_{v_1v_2\cdots v_m}=
Z_{(\partial_iv_1)v_2\cdots v_m}+
Z_{v_1(\partial_iv_2)\cdots v_m}+\cdots +
Z_{v_1v_2\cdots (\partial_iv_m)}.
$$
Obviously, the maps $\delta_i$ are derivations on the algebra ${\cal Z}_0$,
but they are not commuting maps, i.e. the Frobenius property (trivially
verified with the straightforward definition above) $\delta_i\delta_j =
\delta_j\delta_i$, $i\neq j$, is not verified. This comes from the fact
that when one takes the derivatives of an irreducible polynomial $u$,
the polynomials $\partial_i u$, $1\leq i\leq 3$, do not necessarily
factorize out into the same number of factors (a simple example is given
in \cite{DFST}).

A consequence of this fact is the following: If one defines the classical
Nambu bracket on ${\cal Z}_0$ by replacing, in the Jacobian,
the usual product by $\bullet^z$ and the usual partial derivatives
by the maps $\delta_i$, this new bracket will not satisfy the
FI. There will be anomalies in the FI (even at this classical, or
``prequantized" level) due to terms which cannot cancel out each other
because the Frobenius property is not satisfied on ${\cal Z}_0$.

\subsubsection {Taylor series algebra: a solution (case $n =3$)}
Denote by ${\cal E}={\cal Z}_0[y^1,y^2,y^3]$, the algebra of polynomials
in some real variables $(y^1,y^2,y^3)$ with coefficients in ${\cal Z}_0$.
The Taylor series development of translated polynomials $x\mapsto u(x+y)$
can now be written $ u(x+y) = u(x) + \sum_{i}y^i {\partial_i u}(x) +
{1\over 2}\sum_{i,j}y^i y^j {\partial_{ij} u}(x) + \cdots$.
These translated polynomials are multiplied by $(uv)(x+y)=u(x+y)v(x+y)$.
In accordance with our setup we shall look instead at ``Taylor series''
in ${\cal E}$, for $u \in N_1$:
\begin{equation}\label{yyy}
  J(Z_u)=Z_u + \sum_{i}y^i Z_{\partial_i u} + {1\over 2}\sum_{i,j}
y^i y^j Z_{\partial_{ij} u} + \cdots
= \sum_{n} {1 \over {n!}} (\sum_{i}y^i\partial_i)^n(Z_u),
\end{equation}
where $\partial_{i} u$, $\partial_{ij} u$, etc. are the usual derivatives
of $u\in N_1\subset N$ with respect to the variables $x^i$, $x^i$ and
$x^j$, etc., ${\partial_i} Z_u \equiv Z_{\partial_i u}$ and, since in
general the derivatives of $u\in N_1$ are in $N$, one has to factor out
the appropriate constants in $Z_{\partial_i u}$, $Z_{\partial_{ij} u}$, etc.
(i.e. $Z_{\lambda u}\equiv \lambda Z_u$, $u\in N_1$, $\lambda\in \Bbb R$).
$J$ defines an additive map from ${\cal Z}_0$ to ${\cal E}$
(to say that $J$ is multiplicative is tantamount to the Leibniz property).

That algebra is however too large and for our purpose we need to restrict
to the smaller algebra ${\cal A}_0$, the subalgebra of ${\cal E}$ generated
by elements of the form (\ref{yyy}). We shall denote by $\bullet$ the product
in ${\cal A}_0$ which is naturally induced by the product in ${\cal E}$.
In order to define the (classical) Nambu-Poisson structure on ${\cal A}_0$,
we need a correct definition of the derivative of an element
of ${\cal A}_0$. Notice that the derivative ${\partial_i}u(x+y)$ is
again a Taylor series of the form
${\partial_i}u(x) + \sum_j y^j{\partial_{ij}}u(x)+ \cdots$.
We shall thus define the derivative $\Delta_a$, $1\leq a\leq 3$, of an
element of the form (\ref{yyy}) by the natural extension to ${\cal A}_0$
of the previous straightforward definition, i.e.,
\begin{equation}\label{yyya}
\Delta_a(J(Z_u))=J(Z_{\partial_a u})= Z_{\partial_a u} +
\sum_{i}y^i Z_{\partial_{ai} u} + {1\over 2}\sum_{i,j}
y^i y^j Z_{\partial_{aij} u} + \cdots,
\end{equation}
for $u\in N_1$, $1\leq a\leq 3$. One can look at  definition (\ref{yyya})
of $\Delta_a$ as the restriction, to the subset of elements of the form
$J(Z_u)$, of the formal derivative with respect to $y^a$ in the ring
${\cal E}={\cal Z}_0[y^1,y^2,y^3]$. Since
$\Delta_a(J(Z_u))=J(Z_{\partial_a u})$, we have
$\Delta_a({\cal A}_0)={\cal A}_0$ and we get a family of maps
$\Delta_a\colon {\cal A}_0\rightarrow{\cal A}_0$, $1\leq a \leq 3$,
restriction to ${\cal A}_0$ of the derivations with respect to $y^a$,
$1\leq a \leq 3$, in ${\cal E}$. We summarize the properties of $\Delta_a$
in the following:
\begin{lemma}
The maps $\Delta_a\colon {\cal A}_0\rightarrow{\cal A}_0$, $1\leq a\leq 3$,
defined by Eq.~(\ref{yyya}) constitute a family of commuting
derivations (satisfying the Leibniz rule) of the algebra ${\cal A}_0$.
\end{lemma}

The definition of derivatives on ${\cal A}_0$ leads to the following
natural definition of the classical Nambu bracket on the Abelian algebra
${\cal A}_0$ (it is obvious that this defines a Nambu-Poisson structure
on ${\cal A}_0$).
\begin{definition}
The classical Nambu bracket on ${\cal A}_0$ is the trilinear
map taking values in ${\cal A}_0$ given, $\forall A,B,C\in {\cal A}_0$, by:
\begin{equation}\label{xxx}
(A,B,C)\mapsto [A,B,C]_\bullet \equiv \sum_{\sigma\in S_3}
\epsilon(\sigma) \Delta_{\sigma_1}A\bullet \Delta_{\sigma_2}B\bullet
 \Delta_{\sigma_3}C.
\end{equation}
\end{definition}

Now that we have a classical Nambu-Poisson structure on ${\cal A}_0$,
we shall construct a quantum Nambu-Poisson structure by defining a
generalized Abelian deformation $({\cal A}_\nu,\bullet_\nu)$
of $({\cal A}_0,\bullet)$. The construction is based on the map $\alpha$
introduced in Section~3.3.1 and we shall extend the definition of the product
$\bullet_\nu^z$ defined in Section~3.3.2 to the present setting for the
Nambu-Poisson structure on  ${\cal A}_0$.

Let ${\cal E}[\nu]$ be the algebra of polynomials in $\nu$
with coefficients in ${\cal E}$. We consider the subspace  ${\cal A}_\nu$
of ${\cal E}[\nu]$ consisting of series $\sum_{r\geq0}\nu^r A_r$
for which the  coefficient $A_0$ is in ${\cal A}_0$. Then we define a map
$\bullet_\nu \colon {\cal A}_\nu \times {\cal A}_\nu\rightarrow
{\cal E}[\nu]$ by extending the product $\bullet_\nu^z$
defined previously, with $u,v\in N_1$ (it is sufficient to define it
on ${\cal A}_0$ since $\bullet_\nu^z$ annihilates the non-zero powers
of $\nu$):
\begin{equation}\label{product}
J(Z_u)\bullet_\nu J(Z_v) = Z_u\bullet_\nu^z Z_v
+ \sum_{i} y^i ( Z_{\partial_i u}\bullet_\nu^z Z_v +  Z_u\bullet_\nu^z
Z_{\partial_i v})+\cdots .
\end{equation}
Actually $\bullet_\nu$ defines a product on ${\cal A}_\nu$. Moreover,
for $A=\sum_{r\geq0}\nu^r A_r$ and $B=\sum_{s\geq0}\nu^s B_s$ in
${\cal A}_\nu$, we have $A\bullet_\nu B =A_0\bullet_\nu B_0$ and
the coefficient of $\nu^0$ of the latter is $A_0\bullet B_0$ which is in
${\cal A}_0$ since $A_0,B_0\in {\cal A}_0$. This shows that $\bullet_\nu$
is actually a product on ${\cal A}_\nu$ which is Abelian by definition.
Finally, for $\nu=0$, we have $A\bullet_\nu B|_{\nu=0} = A_0\bullet B_0$.
Thus:
\begin{theorem}
The vector space ${\cal A}_\nu$ endowed with the product
$\bullet_\nu$ is a (generalized) Abelian algebra deformation of the
Abelian algebra $({\cal A}_0,\bullet)$.
\end{theorem}

The derivatives $\Delta_a$, $0\leq a\leq 3$, are naturally extended
to ${\cal A}_\nu$. Every element $A\in {\cal A}_\nu$
can be written as $A=\sum_I y^I A_I$, where $I=(i_1,\ldots,i_n)$ is a
multi-index and $A_I\in {\cal Z}_\nu$. Then we have, for
$A,B\in {\cal A}_\nu$,
$A\bullet_\nu B = \sum_{I,J} y^I y^J A_I\bullet^z_\nu B_J$.
Since $({\cal Z}_\nu,\bullet^z_\nu)$ is an Abelian algebra and
the derivative $\Delta_a$ acts as a formal derivative with respect
to $y^a$ on the product $A\bullet_\nu B$, the usual properties
(linearity, Leibniz, Frobenius) of a derivative are still satisfied
on ${\cal A}_\nu$ and we can define the quantum Nambu bracket on
${\cal A}_\nu$.
\begin{definition}
The quantum Nambu bracket on ${\cal A}_\nu$ is the trilinear
map taking values in ${\cal A}_\nu$ defined,
$\forall A,B,C\in {\cal A}_\nu$, by:
\begin{equation}\label{xxxq}
(A,B,C)\mapsto [A,B,C]_{\bullet_\nu} \equiv \sum_{\sigma\in S_3}
\epsilon(\sigma) \Delta_{\sigma_1}A\bullet_\nu \Delta_{\sigma_2}B
\bullet_\nu \Delta_{\sigma_3}C.
\end{equation}
\end{definition}
It is now straightforward to show:
\begin{theorem}
The quantum Nambu bracket (\ref{xxxq}) endows ${\cal A}_\nu$ with
a Nambu-Poisson structure, a (generalized) deformation of the
classical Nambu structure on ${\cal A}_0$.
\end{theorem}

\subsubsection {General case}
The above procedure	can be easily generalized to ${\Bbb R}^n$, $n\geq 2$.
The only non-straightforward modification to be done appears in the
evaluation map. One has to distinguish the cases $n$ even and $n$ odd. If
$n=2p$, $p\geq 1$, then one replaces the partial Moyal product in
(\ref{ta}) by the usual Moyal product on ${\Bbb R}^{2p}$.  If
$n=2p+1$, $p\geq 1$, one uses the partial Moyal product
$\ast_{1\cdots2p}$ on (e.g.) the hyperplane defined by $x_{2p+1}=0$.
The other definitions and properties are directly generalized to
${\Bbb R}^n$. More generally, on a manifold $M$ endowed with a star product
(e.g. a regular Poisson manifold), the same procedure can be developed,
using the a priori given star product in Formula (\ref{ta}).

The question of spectrum of observables can be treated for Zariski products
\cite{DFST} as we have done for deformation quantization. We shall
not develop this point here but shall briefly discuss the question in
Section~4.3.1 together with the ``first quantized" approach to which we
devote Section~4.

Note that the canonical Nambu-Poisson structure of order 2 on ${\Bbb R}^2$
is the usual Poisson structure; there our procedure gives a quantization
of the Poisson bracket ${\cal P}$ different from Moyal, however not on
$N[\nu]$ but on ${\cal A}_\nu$; this quantization will be somewhat like
in field theory.  The same applies to ${\Bbb R}^{2n}$ by starting with a
sum of Poisson brackets on the various ${\Bbb R}^{2}$.

\section {Generalized deformations of the ``first quantized" type;
equivalence, triviality and spectrality \\ for generalized deformations}
At the beginning of Section~3.3.1, we mentioned another possibility for
getting a generalized deformation, already at the level of the algebra
$N$ of real polynomials over ${\Bbb R}^n$. In \cite{DF} this was developed
in the covariant case, which we shall present below. Let us first present
in an even simpler context how this can be done.

Take ${\Bbb R}^2$, with coordinates $p$ and $q$. If $p$ (resp. $q$)
corresponds to an operator $P$ (resp. $Q$), then antistandard ordering
associates to a monomial $q^jp^k$ the operator $Q^jP^k$ and we define in
this way a star product $\ast_S$ which \cite{FS} is nothing but the product
of symbols of pseudodifferential operators. It is equivalent to $\ast_M$
(Moyal) as a star product. The same can be done with standard (first $P$)
and normal (with $P \pm iQ$) orderings, etc. Now change the product
$\times_\alpha$ in (\ref{ta}) into a ``sun" product $\odot_{MS}$ defined
by $(q^{j_1}p^{k_1})\odot_{MS} (q^{j_2}p^{k_2}) =
q^{j_1+j_2} \ast_M p^{k_1+k_2}$ on monomials. This is obviously an Abelian
product but now it has some ``built in" quantization. We extend it by
linearity to all of $N[[\nu]]$ with the convention that it annihilates
nonzero powers of  $\nu$.
We get in this way a generalized deformation, an associative,
distributive and Abelian product with which we can proceed as before, but
this time at the level of the algebra of polynomials on ${\Bbb R}^2$.
If now we consider some ${\Bbb R}^2 \subset {\Bbb R}^3$, we can quantize
Nambu brackets at this ``quantum mechanical" level.	

More generally (and in rough terms), if $\ast$ is some star product on some
manifold $M$ and $N \subset C^\infty(M)$ is some algebra, ${\cal S}$ the
symmetric algebra over it, given some map $\alpha: N[[\nu]] \rightarrow
{\cal S}$ defined in a way similar to the previous case (e.g. by performing
some ordering operation on basis elements of $N$ and extending to formal
series with the convention that it annihilates nonzero powers of $\nu$)
and an evaluation map $T$ defined by replacing $\otimes$ by $\ast$,
we can introduce a ``sun" product $\aa = T\alpha$.
The question which has to be addressed is how trivial is this type of
quantization, but before we shall present the more sophisticated form
studied in \cite{DF}.

\subsection {Covariant star products and Zariski quantization}
Let $N$ be the algebra of real polynomials over ${\Bbb R}^n$ ($n=3$ in the
following) and define a sun product $\aa =T\alpha$ where $\alpha$:
$N[\nu] \rightarrow {\cal S}$ maps the monomial $x_1^{k_1}\cdots x_n^{k_n}
\in N$	into a similar expression belonging to the symmetric tensor algebra
${\cal S}$ over $N$ obtained by replacing the product of the coordinates
$x_i$ ($i=1,\ldots, n$)  by the symmetrized tensor product $\otimes$;
we extend by linearity to $N$ and to $N[\nu]$ by the convention that
$\alpha$ annihilates nonzero powers of $\nu$ (in other words, we extend
by linearity to $N[\nu]$ and compose with the canonical projection of
$N[\nu]$ on $N$). The evaluation map $T$: ${\cal S} \rightarrow N[\nu]$
is defined, for $F_i\in N$, $1\leq i\leq k$, $\forall k\geq 1$, by
\begin{equation}\label{defT}
T(F_1\otimes\cdots\otimes F_k)={1\over k!}\sum_{\sigma\in S_k}
F_{\sigma_1}\ast\cdots\ast F_{\sigma_k}
\end{equation}
where $\ast$ is some star product on $N$. We call $\aa$ the {\it sun
product associated with the star product} $\ast$. It is obviously
Abelian and associative.

In the following we consider ${\Bbb R}^3$ as the dual of the
$\frak{su}(2)$-Lie algebra and take $\ast^{}_M$ to be the Moyal product
on ${\Bbb R}^6$. The Lie algebra ${\frak{su}}(2)$ can be realized with the
functions $L_i(p,q) = \sum_{1\leq j,k \leq 3} \varepsilon_{ijk} p_j q_k$,
$1\leq i,j,k\leq 3$, on ${\Bbb R}^6$, where $\varepsilon_{ijk}$
is the totally skew-symmetric tensor with $\varepsilon_{123}=1$:
\begin{equation}\label{lsu}
 [L_i,L_j]_M^{}\equiv  {1 \over2\nu} (L_i\ast^{}_M L_j  - L_j \ast^{}_M L_i)
=  \sum_{1\leq k \leq 3}\varepsilon_{ijk} L_k.
\end{equation}
One can easily show that for any real polynomial $F \in N({\Bbb R}^3)$,
the polynomial on ${\Bbb R}^6$ defined by $F(L_1,L_2,L_3)$
satisfies $(1\leq i,j,k \leq 3)$
\begin{equation}\label{lproduct}
L_i\ast^{}_M F= L_i F +\nu\sum_{1\leq j,k \leq 3}
\varepsilon_{ijk} L_k {\partial F \over\partial L_j }
+ \nu^2 \Bigl( 2 {\partial F \over\partial L_i} +\sum_{1\leq j \leq 3}
L_j {\partial^2 F \over\partial L_i \partial L_j }\Bigr).
\end{equation}
Therefore $L_{i_1}\ast^{}_M\cdots\ast^{}_M L_{i_k}$ is a polynomial in
$(L_1,L_2,L_3)$ and any polynomial in $(L_1,L_2,L_3)$ can be
expressed as a $\ast^{}_M$-polynomial of the $L_i$'s, so that the product
$F \ast^{}_M G$ of two polynomials $F,G$ in $(L_1,L_2,L_3)$ is again
a polynomial in $(L_1,L_2,L_3)$. Hence from the Moyal product
$\ast^{}_M$ on ${\Bbb R}^6$, we get a star product on ${\Bbb R}^3$
satisfying Eq.~(\ref{lproduct}) for any polynomial $F$. It
is actually an invariant (and covariant) star product on
$\frak{su}(2)^* \sim {\Bbb R}^3$. Herebelow we shall denote by $\ast$
this star product and by $\aa$ the associated $\aa$-product. Since
$L_i\ast L_j = L_i L_j +\nu\sum_{1\leq k \leq 3}\varepsilon_{ijk} L_k
+ 2 \nu^2 \delta_{ij}$, $1\leq i,j \leq 3$,
and $L_i\aa L_j = L_i L_j + 2 \nu^2 \delta_{ij}$, this $\ast$-product
does provide quantum terms.

Let us write  $F\ast G = \sum_{r\geq 0} \nu^r C_r(F,G)$, where
$C_1(F,G)= {\cal P}(F,G)\equiv \sum_{1\leq i,j,k \leq 3}\varepsilon_{ijk}
L_k{\partial F \over\partial L_i}{\partial G \over\partial L_j}$
is the standard Poisson bracket on $\frak{su}(2)^* \sim {\Bbb R}^3$.
We denote by $\Delta$ the Laplacian operator:
$\Delta=\sum_{1\leq k\leq3} {\partial^2\over\partial L_k^2}$.
{}From (\ref{lproduct}) we get
$C_2(L_{i},F) = (2+ {\cal D})({\partial F\over\partial L_{i}})$,
in terms of the differential operator (constant on the subspace
$H_n \subset N$ of homogeneous polynomials of degree $n$)
 ${\cal D} =\sum_{1\leq k \leq 3} L_k{\partial\over\partial L_k}$.
Complicated calculations give $\forall F,G\in	H_n$,
$F\aa G =\sum_{r\geq0} \nu^{2r} a(2n,r) \Delta^r(FG)$,
with $\forall n\geq1$,  $a(n,r) = {1\over n}((n-2r) a(n-1,r) +
(n-2r+2)a(n-1,r-1))$ and $a(n,0) = 1 = a(0,r)$.
After some more calculations \cite{DF} we obtain
\begin{equation}\label{euler}
a(n,r) = \sum_{j_1 + \cdots + j_{n-2r+2}=r\atop j_1,\ldots, j_{n-2r+2}\geq0}
\gamma^{}_{j_1}\gamma^{}_{j_2}\tau^{}_{j_3}\cdots
\tau^{}_{j_{n-2r+2}}, \quad n\geq 2r,
\end{equation}
where  $\gamma^{}_{n}= E_{2n}/(2n)!$ and
$\tau^{}_n = 2^{2n+2}(2^{2n+2}-1) B_{2n+2}/(2n+2)!$, the constants
$E_k$ (resp. $B_k$) being the Euler (resp. Bernoulli) numbers.
Let $B(k,p)$ be the set of partitions of $k$ of length $p$, i.e.
the number of ways to write $k$ as a sum of $p$ strictly positive integers:
$k=n_1+\cdots+n_p$, with $n_1\geq \cdots \geq n_p$, and
$A_k=\sum_{i+j=k\atop i,j \geq0} \gamma^{}_{i}\gamma^{}_{j}$, $k\geq0$.
In the above notations we then get:
\begin{theorem}
{\it The sun product $\aa$ associated with the invariant star pro\-duct on
$\frak{su}(2)^*$ defined by Eq.~(\ref{lproduct}) admits the following form:
\begin{equation}\label{sun}
F\aa G = FG + \sum_{r\geq1} \nu^{2r}\eta_r(FG),
\quad F,G\in N \quad ({polynomials \;\, on}\;\,  {\Bbb R}^3),
\end{equation}
where $\eta_r\colon  N\rightarrow N$, $r\geq1$, are differential
operators given by
\begin{equation}\label{nr}
\eta_r(F)= \left(A_r + \sum_{p=1}^{r} z_{p,r} {\cal D}({\cal D}-1)
\cdots({\cal D}-p+1)\right)
\Delta^r(F), \quad F\in N,
\end{equation}
with $z_{p,r} = \sum_{k=p}^r A_{r-k}\sum_{(n_1,\ldots,n_p)\in B(k,p)}
(m_1!\cdots m_r!)^{-1}\tau^{}_{n_1}\cdots\tau^{}_{n_{p}}$.\\
This $\aa$-product has a unique extension to
$F,G \in C^{\infty}({\Bbb R}^3)$.}
\end{theorem}

\subsection {Two notions of equivalence and triviality}
\subsubsection {Sun products}
 \begin{definition}
Two products $\aa$ and $\aa'$ are said to be A-equivalent,
if there exists a ${\Bbb R}[[\nu]]$-linear (formally invertible) map
$S_\nu\colon N[[\nu]]\mapsto N[[\nu]]$ of the form
$S_\nu=\sum_{r\geq 0} \nu^r S_r$, where the $S_r\colon N\mapsto N$,
$r\geq1$, are differential operators and $S_0=Id$, such that
\begin{equation}\label{Aeq}
 \left. S_\nu(F \aa G) = S_\mu (F) \aa' S_\mu (G)\right|_{\mu= \nu},
\quad \forall F,G \in N.
\end{equation}
We say that they are B-equivalent if the map $S_\nu$ satisfies
\begin{equation}\label{Beq}
S_\nu(F\aa G)= S_\nu(F)\aa' S_\nu(G), \quad \forall F,G \in N.
\end{equation}
\end{definition}
It is straightforward to check that both notions are indeed equivalence
relations (between $\aa$-products). Note that there is no order
relation between both; in particular two $\aa$-products which are
B-equivalent are not necessarily A-equivalent. Relation (\ref{Beq})
can be written as $S_\nu(F\aa G)= F \aa' G$ and relation (\ref{Aeq})
amounts to
\begin{equation}\label{cAeq}
\sum_{r,s\geq0} \nu^{r+s}S_s(\rho_r(FG))
= \sum_{r,s,s'\geq0}\nu^{r+s+s'}\rho_r'(S_s(F)S_{s'}(G)),\quad F,G \in N,
\end{equation}
where $\rho_r$ (resp. $\rho_r'$) are the cochains of the product
$\aa$ (resp. $\aa'$).

Triviality is equivalence with the original product of the algebra, the
usual product. That is:
\begin{definition}
A sun product $\aa$ is said strongly (resp., weakly) trivial if it is
A-equivalent (resp., B-equivalent) with the usual product.
\end{definition}
In other words, if we take for $\aa'$ the usual product,  $\aa$ is
strongly trivial if there exists a (formally invertible)
${\Bbb R}[[\nu]]$-linear map $S_\nu\colon N[[\nu]]\mapsto N[[\nu]]$ of the
form $S_\nu=\sum_{r\geq 0} \nu^r S_r$, where the $S_r\colon N\mapsto N$,
$r\geq1$, are differential operators and $S_0=Id$, such that
(denoting by $\cdot$ the usual product):
\begin{equation}\label{str}
S_\nu(F\aa G)= S_\nu(F)\cdot S_\nu(G), \quad \forall F,G \in N.
\end{equation}
$\aa$ is weakly trivial if (with the same notations)
\begin{equation}\label{wtr}
S_\nu(F\cdot G)	= F\aa G, \quad \forall F,G \in N,
\end{equation}
or equivalently  $F\cdot G = S_\nu^{-1}(F\aa G)$, which corresponds to
taking for $\aa'$ (instead of $\aa$) in (\ref{Beq}) the usual product. We
have here manifest symmetry in the substitution of the usual product.
If in (\ref{Aeq}) we take for $\aa$ the usual product we get
$\left. S_\nu(F \cdot G) = S_\mu (F) \aa S_\mu (G)\right|_{\mu= \nu}$,
which (using (\ref{cAeq})) is equivalent to (\ref{str}) with $S^{-1}$
instead of $S$. From Definition~7 we see that a strongly trivial
sun product is weakly trivial (hence the terminology).
As a matter of fact one proves easily:
\begin{proposition}
A $\aa$-product is strongly trivial if and only if it coincides with the
usual product. A $\aa$-product is weakly trivial whenever its cochains are
given by differential operators.
\end{proposition}
The proof of the first part is by elementary calculation. For the second
part one writes (for $F,G \in N$) $F\aa G = \rho(F\cdot G)$ where
$\rho=Id + \sum_{r\geq1} \nu^r \rho_r$ is formally invertible and the
$\rho_r$'s are differential operators acting on $N$.
{} From Theorem~4 we get:
\begin{corollary}
The sun product constructed in Section~4.1 is weakly trivial but not
strongly trivial.
\end{corollary}
\subsubsection {Zariski products}
The previous treatment (with the obvious modifications) applies to the
case of Zariski quantization. But in this case we see easily that
{\it the Zariski product $\bullet_\nu$ is never trivial in either sense}
(except when the defining product $\ast_\nu$ is not a star product but the
usual product).

It is true that the map $T\alpha$ intertwines the Zariski product
$\bullet_\nu$ with the usual product $\bullet$ on the classical algebra
${\cal A}_0$. But it is neither invertible as a formal series nor acting by
differential operators. Even if one defines $S_\nu$ by taking the
${\Bbb R}[[\nu]]$-linear extension of the restriction of $T\circ\alpha$
to ${\cal A}_0$, one gets an invertible map but not an intertwining map
(in the sense of B-equivalence) as $S_\nu$ is not given by differential
operators on ${\cal A}_0$. We restrict the intertwining operators to be
defined by {\it differentiable} cochains (in spite of the fact that in
general the deformations considered are not defined by differentiable
cochains) because eventually the non-triviality for the Abelian
generalized deformations should be linked with the usual (differentiable)
Hochschild or Harrison cohomologies.

\subsection {Spectrality and concluding remarks}
\subsubsection{Exponentials and spectrum}
As we recalled in Section~3.2, the spectrum of an observable $f$ is obtained,
in deformation quantization, via the $\ast$-exponential (\ref{exp}).
If we take a DrG-deformation of the form (\ref{st}) which is trivial
with intertwining operator $T$ (note that a star product is never trivial
because the Poisson bracket is a nontrivial Hochschild 2-cocycle), the
corresponding exponential is the inverse image of the usual exponential
of $Tf$, hence the spectrum is continuous.

For a $\aa$ product (e.g. the case of $SU(2)$) we can define
$\exp_{\scriptscriptstyle {\odot}}(f)$ like in the star case (\ref{exp})
with $\aa$ instead of $\ast$. This exponential coincides with the
star exponential for linear elements like $L_3$ in the case of $SU(2)$
but not for expressions like 
$H_1=\frac{1}{2}(L_1^2+L_2^2+L_3^2)$. In particular the Fourier
decomposition of  $\exp_{\scriptscriptstyle {\odot}}(L_3)$ is discrete
but for $H_1$ the spectra do not coincide. Now if we take a weakly
trivial $\aa$ product, the intertwiner $S$ is not enough to trivialize
the spectrum, in contradistinction with the previous case.

For the Zariski product $\bullet_\nu^z$ we define similarly the
corresponding exponential. Here, in addition to linear elements, the star
and $\bullet_\nu^z$ exponentials coincide for irreducible polynomials,
giving e.g. the usual discrete spectrum for the harmonic oscillator
Hamiltonian.

\subsubsection{Quantized Nambu bracket}
A quantization of the classical Nambu bracket is achieved by
replacing the usual product by the $\aa$-product of Section~4.
Due to the properties of the $\aa$-product, it is easy to see
that actually the quantized Nambu bracket is given by:
$$
[F,G,H]_{\scriptscriptstyle {\odot}_\nu} = T(\alpha(\{F,G,H\})),
\quad F,G,H\in C^\infty({\Bbb R}^3),
$$
where $\{F,G,H\}$ denotes the classical Nambu bracket on ${\Bbb R}^3$,
i.e., the Jacobian. Though the Leibniz rule is
not satisfied for the
$\aa$-product, this quantized Nambu bracket does satisfy the
Fundamental Identity. Indeed only the weaker form of the Leibniz rule
\begin{eqnarray*}
&&F\aa(\frac{\partial}{\partial L_i}(G\aa H) -
G\aa\frac{\partial H}{\partial L_i} -
\frac{\partial G}{\partial L_i}\aa H)\\
&&\quad=T(\alpha(F(\frac{\partial}{\partial L_i}(G H) -
G\frac{\partial H}{\partial L_i} -\frac{\partial G}{\partial L_i}H))) =0,
\end{eqnarray*}
is required to ensure that the Fundamental Identity is verified by the
quantized Nambu bracket.

Finally let us mention that, as for the $\aa$-product case,
the quantized Nambu bracket is weakly trivial, but not strongly trivial.


\end{document}